\documentclass[11pt,a4paper]{article}
\usepackage[margin=27mm]{geometry}
\usepackage[T1]{fontenc}
\usepackage{lmodern,microtype}
\usepackage{amsmath,amssymb,amsthm,mathtools}
\usepackage{booktabs,tabularx}
\usepackage{xurl}
\usepackage[hidelinks]{hyperref}
\usepackage[numbers,sort&compress]{natbib}
\hypersetup{pdftitle={An Explicit All-Range Reduction for the Linnik--Goldbach Problem with Eighteen Powers of Two},pdfsubject={A binary bridge and a weighted reflection criterion},pdfauthor={Zhengwang Luo; Xidan Song; Contributors from the Goldbach Open Conjecture Platform and the Zhihu AI4S Community}}
\newtheorem{theorem}{Theorem}[section]

\newtheorem{lemma}[theorem]{Lemma}

\theoremstyle{definition}

\theoremstyle{remark}

\numberwithin{equation}{section}
\newcommand{\NN}{\mathbb N}
\newcommand{\ZZ}{\mathbb Z}
\newcommand{\PP}{\mathbb P}
\newcommand{\TT}{\mathbb T}
\newcommand{\GB}{\mathcal G}
\newcommand{\RP}[1]{\mathsf R_{#1}^{+}}
\newcommand{\RLE}[1]{\mathsf R_{\leq #1}^{+}}
\newcommand{\Nbr}{H_{18}}
\newcommand{\epsN}{\varepsilon_N}
\DeclareMathOperator{\pop}{pop}

\title{An Explicit All-Range Reduction\\for the Linnik--Goldbach Problem\\[3pt]
\large with Eighteen Powers of Two\\[0.7em]
{\normalsize\itshape To mark the sixtieth anniversary of Chen Jingrun's announcement of his theorem}}
\author{%
\begin{tabular}{c}
Zhengwang Luo\textsuperscript{1}\qquad
Xidan Song\textsuperscript{2,*}\\[0.6em]
\href{https://www.goldbach.tech/}{Contributors from the Goldbach Open Conjecture Platform}\\
\href{https://www.goldbach.tech/}{and the Zhihu AI4S Community}\\[0.15em]
{\small\url{https://www.goldbach.tech/}}\\[0.7em]
{\small\textsuperscript{1}Sichuan University,
\href{mailto:luozhengwang@stu.scu.edu.cn}{\nolinkurl{luozhengwang@stu.scu.edu.cn}}}\\[0.2em]
{\small\textsuperscript{2}Wuhan University,
\href{mailto:liuqingdamafen@whu.edu.cn}{\nolinkurl{liuqingdamafen@whu.edu.cn}}}\\[0.2em]
{\small\textsuperscript{*}Corresponding author}
\end{tabular}}
\date{28 November 2026}
\begin{document}
\maketitle

\begin{abstract}
We present a constructive reduction for representing every even integer as two primes and at most eighteen positive powers of two. A parameterized binary bridge extends a Goldbach verification through $B$ to exact $r$-power representations through $B+(2^{r+1}-2)2^E$, whenever $1\leq r\leq E$ and $2^E\leq B-4$. For $B=4\cdot10^{18}$, $E=61$, and $r=18$, the resulting endpoint is $H_{18}=1208925207928610747318272$. Above this endpoint, we use the compactly supported weight $\sin^2(\pi x/N)$ and a weighted reflection inequality for the representation of an integer as one prime and nine powers of two. We prove an explicit first-moment estimate under RH and reduce positivity of the eighteen-power representation count to a precisely stated second-moment bound. We also give the prime-power removal, polygonal large-value, and matrix comparison lemmas needed to interpret the accompanying numerical framework. Under GRH and the stated high-range moment bound, the published finite Goldbach verification yields exact eighteen-power representations for every even $N\geq40$, and at-most-eighteen representations for every even $N\geq4$.
\end{abstract}


\section{Introduction}

Linnik's approximation to the binary Goldbach problem permits a bounded number of powers of two in addition to two primes \cite{Linnik1953}. Its asymptotic form asks for a fixed number of powers that suffices for every sufficiently large even integer. An explicit result for all even integers requires an overlap between an analytic range and a finite range where representations are known. The binary bridge below provides such an overlap by an elementary construction.

\subsection{Notation and conventions}

We write $\NN=\{1,2,\ldots\}$, $\NN_0=\{0,1,2,\ldots\}$, and $\PP$ for the primes. All logarithms are natural, except $\log_2 x=(\log x)/(\log2)$. We put $e(t)=\exp(2\pi i t)$ and $\TT=\mathbb R/\mathbb Z$, with normalized Lebesgue measure. The generalized Riemann hypothesis, or GRH, refers to the assertion that every nontrivial zero of every primitive Dirichlet $L$-function has real part $1/2$; it includes the Riemann hypothesis for $\zeta(s)$.

For $r\in\NN_0$, the statement $\RP{r}(N)$ means that there are primes $p,q$ and integers $k_1,\ldots,k_r\geq1$ with
\begin{equation}\label{eq:rep}
 N=p+q+\sum_{j=1}^{r}2^{k_j}.
\end{equation}
We write $\RLE{r}(N)$ if $\RP{s}(N)$ holds for some $0\leq s\leq r$. Both the primes and the exponents may repeat. An empty sum is $0$, so $\RP{0}(N)$ is the binary Goldbach statement for $N$. A \emph{positive power of two} always means $2^k$ with $k\geq1$.

For an even integer $B\geq4$, write
\begin{equation}\label{eq:base}
 \GB(B):\quad
 \text{every even }n\in[4,B]\text{ satisfies }\RP{0}(n).
\end{equation}
The constants used in the eighteen-power specialization are
\begin{equation}\label{eq:constants}
 \begin{gathered}
 B=4\cdot10^{18},\qquad P=2^{61},\\
 \Nbr=B+(2^{19}-2)P=1208925207928610747318272.
 \end{gathered}
\end{equation}

For the related $1+2$ problem, let $\Omega(m)$ count prime factors with multiplicity. We use the squarefree convention
\begin{equation}\label{eq:p2}
 \mathcal P_2^{\mathrm{sf}}
 =\PP\cup\{qr:q,r\in\PP,\ q<r\}.
\end{equation}
Its elements have at most two prime factors; $1$ and prime squares are excluded.

\subsection{Related work and provenance}

Heath-Brown and Puchta obtained thirteen powers unconditionally and seven under GRH for sufficiently large even integers \cite{HBP2002}. Pintz and Ruzsa developed the associated circle-method, singular-series, and lacunary-sum estimates \cite{PR2003,PR2020}. Johnston and Trudgian obtained six powers under GRH \cite{JT2026}; Li and Liu give six unconditionally and five under GRH \cite{LL2026}. These are statements about sufficiently large even integers. The present argument concerns an explicit finite-to-global reduction, in which the crossover is specified as an integer.

Chen proved that every sufficiently large even integer has a representation $N=p+m$ with $p\in\PP$ and $1\leq\Omega(m)\leq2$ \cite[Theorem~I]{Chen1973}. The global $1+2$ result cited here is taken from Li and Shen's \emph{Our Proof} manuscript and repository \cite{LiShenGlobal12}, Theorems~4 and~4.2. Under Dirichlet GRH, their result is
\begin{equation}\label{eq:global12}
 \forall N\geq4\text{ even},\qquad
 N=p+m\quad\text{with }p\in\PP,\ m\in\mathcal P_2^{\mathrm{sf}}.
\end{equation}
Their all-range argument uses the published verification of binary Goldbach through $4\cdot10^{18}$ by Oliveira e Silva, Herzog, and Pardi \cite{OSHP2014}. The same finite statement $\GB(B)$ is the starting point of our binary bridge.

The eighteen-power construction and numerical parameters are those of Visio-Vanitas's manuscript and source snapshot \cite{SourceSnapshot}. We use the branch \path{feat/one-plus-one-v436} at revision \path{a477ae93}; the complete revision is recorded in the bibliography. This includes the code of its parent revision \path{d5081e7a} and supersedes the twenty-two-power baseline submitted in workbench PR~4727 \cite{Workbench4727}. Our exposition gives the general arithmetic bridge, the weighted reflection argument, and explicit interfaces for the analytic estimates.

\subsection{The main statements}

\begin{theorem}[The eighteen-power bridge]\label{thm:bridge18}
Assume $\GB(B)$, with $B$ as in \eqref{eq:constants}. Every even integer $40\leq N\leq\Nbr$ satisfies $\RP{18}(N)$. Every even integer $4\leq N\leq\Nbr$ satisfies $\RLE{18}(N)$.
\end{theorem}

We next specify the analytic input for the remaining range. For $N>\Nbr$, put
\begin{equation}\label{eq:weight}
 \begin{aligned}
 L&=L_N=\left\lfloor\log_2(N/8)\right\rfloor,\qquad c=\frac18,\\
 w_N(x)&=\sin^2(\pi x/N)\mathbf1_{[0,N]}(x),\qquad
 E_2=E_2(N)=\int_{\mathbb R}w_N(x)^2\,dx=\frac{3N}{8}.
 \end{aligned}
\end{equation}
Thus $L\geq76$ and $2^L\leq cN$. Define the finite exponential sums
\begin{equation}\label{eq:sums}
 T_N(\alpha)=\sum_{p\in\PP}(\log p)w_N(p)e(p\alpha),
 \qquad G_L(\alpha)=\sum_{j=1}^{L}e(2^j\alpha),
\end{equation}
and the normalized second moment
\begin{equation}\label{eq:moment-def}
 m_9(N)=\frac1{2E_2L^{18}}
          \int_{\TT}|T_N(\alpha)G_L(\alpha)^9|^2\,d\alpha.
\end{equation}
The explicit first-moment error used below is
\begin{equation}\label{eq:epsilon}
 \epsN=\frac{2\pi^2}{3}
       \left(\frac{3}{16L}+\frac{9}{2L^2}\right)
       +\frac{2(\log N)^2}{\pi\sqrt N}.
\end{equation}

\begin{theorem}[All-range reduction]\label{thm:global}
Assume GRH, $\GB(B)$, and the high-range moment bound
\begin{equation}\label{eq:moment-input}
 m_9(N)<2(1-\epsN)^2
 \qquad\text{for every even integer }N>\Nbr.
\end{equation}
Then every even integer $N\geq40$ satisfies $\RP{18}(N)$, and every even integer $N\geq4$ satisfies $\RLE{18}(N)$.
\end{theorem}

The endpoint $40$ is necessary in the exact statement, since the smallest possible sum of two primes and eighteen positive powers of two is $2+2+18\cdot2=40$. Theorem~\ref{thm:global} isolates a second-moment estimate, while the first-moment estimate is proved in Section~\ref{sec:reflection}. Section~\ref{sec:analytic} describes the analytic and numerical structure proposed for \eqref{eq:moment-input}.

\section{A parameterized binary bridge}\label{sec:bridge}

For $M\in\NN_0$, let $\pop(M)$ be the number of nonzero digits in its binary expansion.

\begin{lemma}\label{lem:bits}
If $r\geq1$ and $1\leq M\leq2^{r+1}-2$, then $1\leq\pop(M)\leq r$.
\end{lemma}

\begin{proof}
The expansion has at most $r+1$ binary positions. The only integer below $2^{r+1}$ with all these positions nonzero is $2^{r+1}-1$, which is excluded. At least one position is nonzero because $M>0$.
\end{proof}

\begin{lemma}[Splitting powers]\label{lem:split}
Suppose $1\leq s\leq r\leq E$ and $U=\sum_{j=1}^{s}2^{a_j}$ with $a_j\geq E$. Then $U$ is a sum of exactly $r$ positive powers of two, with every final exponent at least $E-r+s$.
\end{lemma}

\begin{proof}
Replace a term $2^a$ by $2^{a-1}+2^{a-1}$, repeating this operation $r-s$ times. Each replacement preserves the sum and adds one term. After $t$ replacements, every exponent is at least $E-t$. Before a further replacement is needed, $t<r-s$, so $E-t\geq E-r+s+1\geq2$. Every step is therefore valid with positive final exponents.
\end{proof}

\begin{theorem}[General bridge]\label{thm:general-bridge}
Let $B$ be even and let $r,E$ be integers satisfying
\[
 1\leq r\leq E,\qquad 2^E\leq B-4.
\]
Assume $\GB(B)$ and put $H_r(B,E)=B+(2^{r+1}-2)2^E$. Then every even integer
\[
 2r+4\leq N\leq H_r(B,E)
\]
satisfies $\RP{r}(N)$.
\end{theorem}

\begin{proof}
Write $P=2^E$. If $N\leq B+2r$, the even integer $N-2r$ lies in $[4,B]$. Write it as $p+q$ by $\GB(B)$ and add $r$ copies of $2$.

If $N>B+2r$, set
\begin{equation}\label{eq:ceiling}
 A=N-B,\qquad M=\left\lceil\frac AP\right\rceil,
 \qquad n=N-MP.
\end{equation}
The upper bound on $N$ implies $1\leq M\leq2^{r+1}-2$, and the ceiling gives
\[
 0\leq MP-A<P,\qquad B-P<n\leq B.
\]
Since $P\leq B-4$, the even integer $n$ lies in $[4,B]$; choose primes $p,q$ with $n=p+q$.

By Lemma~\ref{lem:bits}, write
\[
 M=\sum_{j=1}^{s}2^{b_j},\qquad
 1\leq s\leq r,\quad b_j\geq0.
\]
Then $MP=\sum_{j=1}^{s}2^{b_j+E}$. Lemma~\ref{lem:split} expresses $MP$ as exactly $r$ positive powers of two, and $N=p+q+MP$ completes the proof.
\end{proof}

\begin{proof}[Proof of Theorem~\ref{thm:bridge18}]
Take $r=18$, $E=61$, and $B=4\cdot10^{18}$. The inequalities in Theorem~\ref{thm:general-bridge} hold because
\[
 18\leq61,\qquad 2^{61}=2305843009213693952\leq B-4.
\]
The resulting upper endpoint is exactly $\Nbr$ in \eqref{eq:constants}. For even $4\leq N<40$, put $s=(N-4)/2\leq17$ and use
\begin{equation}\label{eq:small}
 N=2+2+\underbrace{2+\cdots+2}_{s\text{ terms}}.
\end{equation}
This includes the empty sum when $N=4$.
\end{proof}

For the eighteen-power bridge, no more than seventeen splittings are needed, so all exponents produced in the second case are at least $61-17=44$. The finite construction uses no analytic estimate. Its role is to place the beginning of the remaining range at $L=76$ in the normalization \eqref{eq:weight}.

\section{Weighted reflection and the first moment}\label{sec:reflection}

\subsection{From nine powers to eighteen powers}

Fix an even $N>\Nbr$ and retain \eqref{eq:weight}--\eqref{eq:sums}. For a tuple $\boldsymbol{j}=(j_1,\ldots,j_k)\in\{1,\ldots,L\}^k$, put $\sigma(\boldsymbol{j})=\sum_{i=1}^{k}2^{j_i}$. Define the nonnegative, finitely supported function
\begin{equation}\label{eq:f}
 f_k(n)=\sum_{\substack{p\in\PP,\ \boldsymbol{j}\in\{1,\ldots,L\}^k\\
                       p+\sigma(\boldsymbol{j})=n}}
                    (\log p)w_N(p),\qquad n\in\ZZ.
\end{equation}
Its Fourier polynomial is $T_NG_L^k$. We use
\begin{equation}\label{eq:A-M-R}
 \begin{aligned}
 A_k&=\sum_{\substack{1\leq n<N\\n\text{ odd}}}f_k(n)w_N(n),
 &M_k&=\sum_{n\in\ZZ}f_k(n)^2,\\
 R_{2k}(N)&=\sum_{n\in\ZZ}f_k(n)f_k(N-n).
 \end{aligned}
\end{equation}
Fourier orthogonality gives
\begin{equation}\label{eq:fourier}
 \begin{aligned}
 M_k&=\int_{\TT}|T_NG_L^k|^2\,d\alpha,\\
 R_{2k}(N)&=\int_{\TT}T_N(\alpha)^2G_L(\alpha)^{2k}e(-N\alpha)\,d\alpha\\
 &=\sum_{\substack{p,q\in\PP,\ 1\leq j_1,\ldots,j_{2k}\leq L\\
                    p+q+\sum_i2^{j_i}=N}}
       (\log p)(\log q)w_N(p)w_N(q).
 \end{aligned}
\end{equation}
The last sum is over ordered tuples, and every term is positive. Thus $R_{2k}(N)>0$ gives $\RP{2k}(N)$ with all exponents at most $L$. The second moment $M_k$ and the representation coefficient $R_{2k}$ are connected by the following weighted version of the reflection argument of \cite[Lemma~2.1]{LL2026}.

\begin{lemma}[Weighted reflection]\label{lem:reflection}
Let $N\geq6$ be even and let $L,k\geq1$ be integers. With the quantities defined in \eqref{eq:sums}, \eqref{eq:f}, and \eqref{eq:A-M-R},
\begin{equation}\label{eq:reflection}
 R_{2k}(N)\geq\frac{4A_k^2}{E_2}-M_k.
\end{equation}
Consequently, with $a_k=A_k/(E_2L^k)$, $m_k=M_k/(2E_2L^{2k})$, and $J_{2k}=R_{2k}/(2E_2L^{2k})$,
\begin{equation}\label{eq:normalized-reflection}
 J_{2k}\geq2a_k^2-m_k.
\end{equation}
\end{lemma}

\begin{proof}
Let $U=\{n:1\leq n<N,\ n\text{ odd}\}$. The identity
\[
 \sin^4 t=\frac{3-4\cos(2t)+\cos(4t)}8
\]
and finite geometric sums give
\begin{equation}\label{eq:odd-energy}
 \sum_{n\in U}w_N(n)^2=\frac{3N}{16}=\frac{E_2}{2}.
\end{equation}
Here $N/2\geq3$, so both nonconstant cosine sums vanish. Also $w_N(N-n)=w_N(n)$, and reflection permutes $U$. Hence
\[
 2A_k=\sum_{n\in U}\bigl(f_k(n)+f_k(N-n)\bigr)w_N(n).
\]
Cauchy--Schwarz and \eqref{eq:odd-energy} imply
\[
 \sum_{n\in U}\bigl(f_k(n)+f_k(N-n)\bigr)^2
 \geq\frac{8A_k^2}{E_2}.
\]
Expanding the left side and using reflection once more gives
\[
 \sum_{n\in U}f_k(n)f_k(N-n)
 \geq\frac{4A_k^2}{E_2}-\sum_{n\in U}f_k(n)^2.
\]
The full representation sum is at least its restriction to $U$, and the last sum is at most $M_k$. This proves both assertions.
\end{proof}

\subsection{The effect of the binary shifts}

\begin{lemma}[Translation estimate]\label{lem:translation}
For every real $h$,
\begin{equation}\label{eq:translation}
 \int_{\mathbb R}w_N(x)w_N(x+h)\,dx
 \geq E_2\left(1-\frac{2\pi^2h^2}{3N^2}\right).
\end{equation}
If $\boldsymbol{j}$ is uniform on $\{1,\ldots,L\}^k$, then
\begin{equation}\label{eq:shift-moment}
 \frac{\mathbb E\sigma(\boldsymbol{j})^2}{N^2}
 \leq\frac{4kc^2}{3L}+\frac{4k(k-1)c^2}{L^2}.
\end{equation}
\end{lemma}

\begin{proof}
Extend $w_N$ by zero outside $[0,N]$. It is continuously differentiable and
\[
 \|w_N'\|_2^2=\frac{\pi^2}{2N},\qquad
 \|w_N(\cdot+h)-w_N\|_2\leq |h|\|w_N'\|_2.
\]
The latter inequality follows by integrating the derivative and applying the triangle inequality in $L^2$. Since
\[
 \int w_N(x)w_N(x+h)\,dx
 =E_2-\frac12\|w_N(\cdot+h)-w_N\|_2^2,
\]
we obtain \eqref{eq:translation}.

For a uniform $j\in\{1,\ldots,L\}$, geometric summation and $2^L\leq cN$ give
\[
 \mathbb E2^j\leq\frac{2cN}{L},\qquad
 \mathbb E4^j\leq\frac{4c^2N^2}{3L}.
\]
Expanding the square of a sum of $k$ independent copies proves \eqref{eq:shift-moment}.
\end{proof}

The translation estimate holds even when $h\geq N$, in which case the correlation integral vanishes and the lower bound may be negative. Thus no restriction on an individual nine-power shift is needed when averaging \eqref{eq:translation}.

\subsection{An explicit first-moment estimate}

Write $\vartheta(x)=\sum_{p\leq x}\log p$ and $\vartheta_{\mathrm{odd}}(x)=\sum_{2<p\leq x}\log p$.

\begin{lemma}[First moment under RH]\label{lem:first-moment}
Assume RH. For every even $N>\Nbr$,
\begin{equation}\label{eq:first-moment}
 \frac{A_9}{E_2L^9}\geq1-\epsN,
 \qquad 0<\epsN<\frac1{40},
\end{equation}
where $\epsN$ is defined in \eqref{eq:epsilon}.
\end{lemma}

\begin{proof}
Schoenfeld's explicit estimate \cite[Theorem~10, (6.3)]{Schoenfeld1976} states that RH implies
\begin{equation}\label{eq:schoenfeld}
 |\vartheta(x)-x|<\frac{\sqrt{x}(\log x)^2}{8\pi}
 \qquad(x\geq599).
\end{equation}
Put $D_N=\sqrt N(\log N)^2/(8\pi)$. On $0\leq x<599$, the elementary estimate $\vartheta(x)\leq x\log x$ for $x\geq2$, and direct treatment of $x<2$, give $|\vartheta(x)-x|<4000$. Since $N>\Nbr>10^{24}$, we have $D_N>4000$. Monotonicity of $\sqrt{x}(\log x)^2$ on $x\geq599$ therefore gives
\[
 \sup_{0\leq x\leq N}|\vartheta_{\mathrm{odd}}(x)-x|
 \leq D_N+\log2\leq2D_N.
\]

For $h\geq0$, let $g_h(x)=w_N(x)w_N(x+h)$. If $h<N$, this function is log-concave in the interior of its support, hence unimodal, with $0\leq g_h\leq1$ and zero endpoint values. Thus its total variation is at most $2$. For $h\geq N$ it vanishes. Partial summation gives, uniformly for all $h\geq0$,
\begin{equation}\label{eq:weighted-pnt}
 \left|\sum_{p>2}(\log p)g_h(p)-\int_{\mathbb R}g_h(x)\,dx\right|
 \leq4D_N
 \leq E_2\frac{2(\log N)^2}{\pi\sqrt N}.
\end{equation}

Each binary shift is even, so the odd indices in $A_9$ correspond exactly to odd primes. Apply \eqref{eq:weighted-pnt} with $h=\sigma(\boldsymbol{j})$, average over all $L^9$ tuples, and use Lemma~\ref{lem:translation}. With $k=9$ and $c=1/8$, \eqref{eq:shift-moment} becomes
\[
 \frac{\mathbb E\sigma^2}{N^2}
 \leq\frac3{16L}+\frac9{2L^2},
\]
which proves the first assertion.

For the numerical bound, $L\geq76$ and $\pi^2<10$ bound the shift term by
\[
 \frac{20}{3}\left(\frac3{16\cdot76}
                   +\frac9{2\cdot76^2}\right)<\frac{11}{500}.
\]
The remaining term is less than $10^{-8}$ for $N>\Nbr$, since $(\log x)^2/\sqrt x$ decreases for $x>\exp(4)$. Their sum is less than $1/40$.
\end{proof}

\begin{proof}[Proof of Theorem~\ref{thm:global}]
For even $N>\Nbr$, Lemmas~\ref{lem:reflection} and~\ref{lem:first-moment}, together with \eqref{eq:moment-input}, imply
\[
 \frac{R_{18}(N)}{2E_2L^{18}}
 \geq2(1-\epsN)^2-m_9(N)>0.
\]
The nonnegative coefficient formula \eqref{eq:fourier} gives $\RP{18}(N)$. Theorem~\ref{thm:bridge18} covers $40\leq N\leq\Nbr$, and \eqref{eq:small} covers the remaining even integers $4\leq N<40$ with at most seventeen powers.
\end{proof}

\section{The analytic second-moment framework}\label{sec:analytic}

The v2 computation \cite{SourceSnapshot} is organized around a normalized bound
\begin{equation}\label{eq:budget-goal}
 m_9(N)\leq1+\mathcal B_N,
 \qquad 4\epsN+\mathcal B_N<1.
\end{equation}
These two inequalities imply \eqref{eq:moment-input}, since
$1+\mathcal B_N<2-4\epsN\leq2(1-\epsN)^2$.
The factor $4\epsN$ comes directly from the reflection argument. We give the exact identities behind the remaining terms.

\subsection{Prime pairs and the use of the sieve}

Put $z(\alpha)=G_L(\alpha)/L$. The monomial case of the polynomial framework is
\[
 P(z)=z^9,\qquad Q(z)=z^8,\qquad P(z)=zQ(z).
\]
The representation coefficient uses $G_L^{18}$, while the second moment uses $|G_L|^{18}$. All powers in the representation therefore occur with positive sign.

The two-power auxiliary moment has the exact expansion
\begin{equation}\label{eq:pair}
 \int_{\TT}|T_NG_L|^2\,d\alpha
 =\sum_{1\leq j,l\leq L}
   \sum_{\substack{p,q\in\PP\\p-q=2^l-2^j}}
   (\log p)(\log q)w_N(p)w_N(q).
\end{equation}
The diagonal is $L\sum_p(\log p)^2w_N(p)^2$. The off-diagonal terms are weighted prime-pair counts, to which the Selberg upper-bound sieve applies \cite{Selberg1947,PR2003,PR2020}. The finite Mellin estimates for this particular weight are specified in \cite{SourceSnapshot}.

The averaged singular-series constant used there is bounded by $1.279$. In its original normalization, Pintz's bound \cite[Lemma~$2''$]{Pintz2006}, quoted by Platt and Trudgian \cite[p.~55]{PT2015}, is
\[
 C_0C_2<1.2784421041<1.279.
\]
Here $C_0=\prod_{p>2}(1-(p-1)^{-2})$; the factor called $C_2$ in \cite{PT2015} is called $R_0$ in the v2 code. The separate asymptotic benchmark uses the sieve constant $6.7814$ arising from Lichtman's work \cite[Theorem~1.2]{Lichtman2023}, with the prime-difference formulation discussed in \cite[Proposition~2]{JT2026}. The finite eighteen-power reduction uses its own explicit sieve bounds.

\subsection{Removing higher prime powers}

Let $\Lambda(n)=\log p$ if $n=p^a$ for a prime $p$ and an integer $a\geq1$, and let $\Lambda(n)=0$ otherwise. Set
\[
 S_N(\alpha)=\sum_{n\leq N}\Lambda(n)w_N(n)e(n\alpha),
 \qquad \psi(x)=\sum_{n\leq x}\Lambda(n).
\]
Let $J_{18}^{\PP}=R_{18}/(2E_2L^{18})$, and define $J_{18}^{\Lambda}$ by replacing $T_N^2$ with $S_N^2$ in its Fourier integral. Likewise, let $m_9^{\Lambda}$ be obtained from \eqref{eq:moment-def} by replacing $T_N$ with $S_N$.

\begin{lemma}[Prime-power removal]\label{lem:prime-powers}
For $N>\Nbr$, both differences satisfy
\begin{equation}\label{eq:pp-removal}
 \begin{aligned}
 0\leq J_{18}^{\Lambda}(N)-J_{18}^{\PP}(N)&\leq\delta_N,\\
 0\leq m_9^{\Lambda}(N)-m_9(N)&\leq\delta_N,
 \end{aligned}
 \qquad
 \delta_N=\frac{\log N}{E_2}\bigl(\psi(N)-\vartheta(N)\bigr)
          <2\cdot10^{-10}.
\end{equation}
\end{lemma}

\begin{proof}
Expand either full-circle integral by orthogonality. For each fixed eighteen-index tuple, a term in the difference has at least one integer that is a proper prime power. Since $0\leq w_N\leq1$ and $\Lambda(m)\leq\log N$ for $m\leq N$, the two possible locations contribute at most
\[
 2\log N\sum_{\substack{p^a\leq N\\a\geq2}}\log p
 =2\log N\bigl(\psi(N)-\vartheta(N)\bigr).
\]
There are $L^{18}$ exponent tuples. Division by $2E_2L^{18}$ proves the upper bounds; nonnegativity follows from the coefficient expansions.

Rosser and Schoenfeld \cite[Theorem~13, (3.36)]{RS1962} give
\[
 \psi(x)-\vartheta(x)<1.42620\sqrt x\qquad(x>0).
\]
Consequently, $\delta_N<3.8032(\log N)/\sqrt N$. The last function decreases for $N>\exp(2)$. Using $\Nbr>1.2\cdot10^{24}$, $\log(1.2\cdot10^{24})<56$, and $\sqrt{1.2\cdot10^{24}}>1.09\cdot10^{12}$ gives the stated bound.
\end{proof}

In particular, a full-circle upper bound for $m_9^{\Lambda}$ also bounds $m_9$. For restricted arc integrals, the approximation error is instead controlled before integration by the corresponding $L^2$ estimates.

\subsection{A coupled major- and minor-arc inequality}

Let $\mathfrak M\subset\TT$ be a measurable major-arc set and $\mathfrak m=\TT\setminus\mathfrak M$. Define the positive measure
\[
 d\nu_N(\alpha)=\frac{|T_N(\alpha)|^2}{2E_2}\,d\alpha,
 \qquad C_{\mathrm{pair}}(N)=\int_{\TT}|z(\alpha)|^2\,d\nu_N(\alpha).
\]
For $0\leq h\leq1$, put
\[
 \mathcal E_h=\{\alpha\in\mathfrak m:|z(\alpha)|^{16}>h\}.
\]

\begin{lemma}\label{lem:coupled}
For every $h\in[0,1]$,
\begin{equation}\label{eq:coupled}
 \begin{split}
 m_9(N)\leq{}&hC_{\mathrm{pair}}(N)
 +\int_{\mathfrak M}\bigl(|z|^{18}-h|z|^2\bigr)\,d\nu_N\\
 &+\int_{\mathcal E_h}|z|^{18}\,d\nu_N.
 \end{split}
\end{equation}
\end{lemma}

\begin{proof}
Subtract $h|z|^2$ from $|z|^{18}$ on the whole circle. The difference is nonpositive on $\mathfrak m\setminus\mathcal E_h$, and is at most $|z|^{18}$ on $\mathcal E_h$. Integration proves the claim.
\end{proof}

The subtraction on the major arcs is coupled to the full prime-pair moment. Retaining it in the form \eqref{eq:coupled} specifies which singular-series terms must be estimated together. At rational points, the finite values $G_L(a/q)/L$ lead to exact character and Ramanujan-sum tables. The v2 computation uses the forty-three odd squarefree moduli $3\leq q\leq105$, followed by a truncation at $2000$ for the remaining arithmetic tail \cite{SourceSnapshot}.

The zero-sum estimates needed for major-arc approximation use the explicit Dirichlet zero counts of Bennett, Martin, O'Bryant, and Rechnitzer \cite[Theorem~1.1]{BMOR2021}. Their theorem applies to primitive characters of conductor $q>1$ at heights $T\geq5/7$ and counts both signs of the zero ordinate. For conductor one, the relevant zeta estimate is due to Hasanalizade, Shen, and Wong \cite[Corollary~1.2]{HSW2022}, valid for $T\geq e$; its positive-ordinate count is doubled when a two-sided count is required. The related explicit arithmetic-progression bounds are recorded separately in \cite{BMOR2018}. These inputs have different roles from the finite rational singular-series tables.

\subsection{Polygonal large-value bounds and the infinite tail}

The large-value analysis of lacunary sums follows the circle-method developments in \cite{HBP2002,PR2003,PR2020,LL2026}. The v2 implementation uses sixteen rational directions and nested polygons to bound the image of $z(\alpha)$ outside a small exceptional set.

\begin{lemma}[Polygonal maximum]\label{lem:polygon}
Let $K\subset\mathbb C$ be a compact convex polygon with vertex set $V(K)$. For every integer $j\geq1$,
\[
 \max_{z\in K}|z|^{2j}=\max_{v\in V(K)}|v|^{2j}.
\]
\end{lemma}

\begin{proof}
The function $z\mapsto|z|^{2j}$ is convex. Express a point of $K$ as a convex combination of its vertices and apply convexity. The reverse inequality follows because every vertex lies in $K$.
\end{proof}

Thus rational vertices permit exact evaluation of the maxima of $|Q(z)|^2=|z|^{16}$ and $|P(z)|^2=|z|^{18}$. If $K_0\subset\cdots\subset K_s$, the corresponding maxima can be combined with bounds for the exceptional measures by finite summation by parts. Explicitly, suppose an integrand $F$ satisfies $F\leq b_i$ on the preimage of $K_i$, where $b_0\leq\cdots\leq b_{s+1}$ and $F\leq b_{s+1}$ everywhere. Then
\begin{equation}\label{eq:layers}
 F\leq b_0+\sum_{i=0}^{s}(b_{i+1}-b_i)
                \mathbf1_{\{z(\alpha)\notin K_i\}}.
\end{equation}
Integrating \eqref{eq:layers} yields the required finite-layer bound.

For an unbounded range of $L$, the transfer matrices admit a direct nonnegative-vector certificate, the comparison underlying the Collatz--Wielandt method \cite{Wielandt1950}.

\begin{lemma}[Matrix comparison]\label{lem:matrix}
Let $A$ be a matrix with nonnegative entries. If $v$ has strictly positive entries and $Av\leq\rho v$ coordinatewise, then $A^nv\leq\rho^nv$ for every $n\geq0$. More generally, $0\leq u\leq Cv$ implies $A^nu\leq C\rho^nv$.
\end{lemma}

\begin{proof}
Multiplication by a nonnegative matrix preserves coordinatewise inequalities. Both assertions follow by induction.
\end{proof}

In applying this lemma to a tail estimate, the matrix growth and the normalization in the exponential-moment bound must be compared for every remaining $L$. The archived calculation records a resulting growth ratio at most $0.916314$ for its chosen tail parameters.

\subsection{The v2 parameter ranges}

For each integer $L\geq76$, define the relevant interval by
\[
 I_L=\{N>\Nbr:8\cdot2^L\leq N<8\cdot2^{L+1}\}.
\]
The finite computation treats $I_{76},\ldots,I_{119}$; the uniform tail begins at $8\cdot2^{120}$. The first four intervals are each subdivided into two pieces. For $76\leq L\leq84$, the source uses a transfer matrix with the specified dyadic exclusions; from $L=85$ onward it uses the pure transfer matrix. The recorded parameters are summarized in Table~\ref{tab:parameters}.

\begin{table}[htbp]
\centering
\small
\caption{Selected values in the archived v2 calculation \cite{SourceSnapshot}. The last column is its normalized budget, denoted by $\rho_N$ in \eqref{eq:source-budget}.}
\label{tab:parameters}
\begin{tabular}{@{}lll@{}}
\toprule
Range & Sieve parameter or bound & Recorded upper budget\\
\midrule
$L=76$ & $z_{\rm sieve}=28000$, $C_{\rm Sel}\leq5.025$ & $\rho_N\leq0.9883$\\
$L=85$ & $C_{\rm Sel}\leq5.383$ & $\rho_N\leq0.9509$\\
$L=119$ & $C_{\rm Sel}\leq7.411$ & $\rho_N\leq0.7235$\\
$L\geq120$ & $C_{\rm Sel}\leq9.0945$ & $\rho_N\leq0.8886$\\
\bottomrule
\end{tabular}
\end{table}

In the source notation, the budget has the structure
\begin{equation}\label{eq:source-budget}
 \begin{split}
 \rho_N={}&4\epsN+A_{\rm ex}+(1-h_0)T_{11}\\
 &+h_0\bigl((C_{\rm Sel}-1)C+D_{\rm diag}+\tau\bigr)
       +X+E_{\rm maj}+E_{\rm dyad},\qquad C=1.279.
 \end{split}
\end{equation}
Here $A_{\rm ex}$ is the contribution of the selected arithmetic moduli, $T_{11}$ is the remaining quadratic arithmetic tail, $h_0$ is a polygonal bound for $|Q(z)|^2$, $D_{\rm diag}$ is the diagonal correction, and $\tau$ is the singular-series truncation term. The source assigns the terms $X$, $E_{\rm maj}$, and $E_{\rm dyad}$ to the exceptional-set, major-arc, and additional dyadic-arc contributions. The actual evaluation replaces the single-$N$ quantities by interval envelopes.

The analytic estimate that connects this budget to Theorem~\ref{thm:global} is
\begin{equation}\label{eq:budget-interface}
 m_9(N)\leq1+\rho_N-4\epsN.
\end{equation}
Together with $\rho_N<1$, it is a sufficient instance of \eqref{eq:budget-goal}. In particular, the coupled term in \eqref{eq:coupled} must justify the coefficient of $T_{11}$ in \eqref{eq:source-budget} together with the remaining prime-pair term. Equations \eqref{eq:moment-input} and \eqref{eq:budget-interface} state the high-range input in mathematical, rather than implementation, terms.

\section{Concluding remarks and scope}

The elementary construction extends the finite Goldbach base to exact eighteen-power representations through $\Nbr\approx1.208925\cdot10^{24}$. The weighted reflection inequality then reduces the infinite range to \eqref{eq:moment-input}. The first-moment estimate, the treatment of every binary shift, and the prime-power removal have been proved above with explicit constants. They make the quantitative target for the second moment precise.

The numerical framework in \cite{SourceSnapshot} replaces the earlier Laplace-weight calculation by the compactly supported $\sin^2$ weight, exact arithmetic tables, and rational polygonal bounds. The two former additive allowances, $0.06$ and $4.00$, are absent from the v2 expressions. Re-evaluation of the pinned script reproduces its archived output, including the smallest reported finite-range margin $0.0117$ and the tail margin $0.1114$. These are properties of the recorded numerical calculation. Establishing the uniform analytic inequality \eqref{eq:budget-interface}, including its coupled major-arc subtraction and all interval and tail enclosures, remains the input required to turn the all-range reduction into a theorem under GRH alone together with the finite Goldbach verification. The numerical margins are not substituted for that analytic implication in the proofs above.

The pinned Lean development \cite{Lean4,Mathlib2020,SourceSnapshot} proves the eighteen-power binary-bridge implication from the finite base, and combines it with the high-range assertion. Both representation predicates now impose positive exponents directly. The exported theorem receives its high-range assertion through \path{LinnikCircleAboveBridgePremise} in \path{OnePlusOne/External/Premises.lean}. Thus its present scope is the arithmetic bridge and all-range composition relative to this external analytic input. The general bridge, the first-moment argument, and the moment inequalities stated here provide separate mathematical lemmas for further formalization.

\section*{Acknowledgments}

We thank Zhihu for providing a forum for our discussions and the many members of the online community who encouraged and supported our work.
We are grateful for the assistance of the AI models used in this research and thank the teams and companies that developed them.
We pay special tribute to the late Professor Chen Jingrun, whose pioneering work on the Goldbach conjecture continues to inspire us.

\section*{AI-use declaration}

This article presents natural-language proofs based on the Lean proofs produced using large language models as part of a project undertaken by Contributors from the Goldbach Open Conjecture Platform and the Zhihu AI4S Community. Jiang Zhexin from The Chinese University of Hong Kong, Shenzhen helped complete the writing of this paper.

\begingroup
\renewcommand{\bibfont}{\small}
\interlinepenalty=10000

\endgroup
\end{document}